%% file: zeta-rational8.tex
\documentclass{amsart}
\usepackage{epsfig}
\usepackage{amscd}
\usepackage{amsmath, amssymb, amsthm}
\usepackage{graphics}
\usepackage{color}
\newtheorem{theorem}{Theorem}
\newtheorem{proposition}{Proposition}[section]
\newtheorem{lemma}[proposition]{Lemma}
\newtheorem{corollary}[theorem]{Corollary}
\newtheorem{conjecture}[theorem]{Conjecture}

\theoremstyle{definition}

\newtheorem*{remark}{Remark}

\numberwithin{equation}{section}

\def\ZZ{{\mathbb Z}}
\def\reals{{\mathbb R}}
\def\cx{{\mathbb C}}
\def\LL{{\mathcal L}}
\def\JJ{{\mathcal J}}

\def\Re{\mathrm{Re}\,}
\def\Im{\mathrm{Im}\,}

\def\diam{\mathrm{diam}\,}
\def\tr{\mathrm{tr}\,}

\def\supp{\mathrm{supp}\,}

\def\Ph{\widehat{\varphi}}
\def\Phd{\widehat{\varphi}_{\gamma,d}}
\def\P{\varphi}
\def\Pd{\varphi_{\gamma,d}}
\def\O{{\mathcal O}}
\begin{document}
\title[Zeta Function for Hyperbolic Rational Maps]{Growth and Zeros of the Zeta Function for Hyperbolic Rational Maps}
\author{Hans Christianson}
\address{Department of Mathematics, University of California, Berkeley, CA 94720 USA}
\email{hans@math.berkeley.edu}
\thanks{The author would like to thank Maciej Zworski for copious help and guidance while writing this paper.  He would also like to thank Daniel Tataru, Mike Christ, Jeremy Marzuola, and Richard Burstein for many helpful conversations and suggestions, as well as the University of California and the National Science Foundation for support.  }
\begin{abstract}
This paper describes new results on the growth and zeros of the Ruelle zeta function for the Julia set of a hyperbolic rational map.  It is shown that the zeta function is bounded by $\exp(C_K |s|^{\delta})$ in strips $|\Re s| \leq K$, where $\delta$ is the dimension of the Julia set.  This leads to bounds on the number of zeros in strips (interpreted as the Pollicott-Ruelle resonances of this dynamical system).  An upper bound on the number of zeros in polynomial regions $\{|\Re s | \leq |\Im s|^\alpha\}$ is given, followed by weaker lower bound estimates in strips $\{\Re s > -C, |\Im s|\leq r\}$, and logarithmic neighbourhoods $\{ |\Re s | \leq \rho \log |\Im s| \}$.  Recent numerical work of Strain-Zworski suggests the upper bounds in strips are optimal.
\end{abstract}
\maketitle

\section{Introduction}
\indent The motivation for the estimates described in this paper comes from scattering 
resonances.  In the case where the underlying fractal set is the limit set of a convex co-compact Schottky group there is a correspondence between zeros of the zeta function and scattering resonances of the classically trapped set (see \cite{GLZ, Ruelle1, Zworski1, Zworski2} for details).  \\
\indent In the case of the Julia set, we are primarily interested in counting zeros of the zeta function $Z$, which may be interpreted as the Pollicott-Ruelle resonances of this dynamical system.  The most interesting case is the number of zeros in regions $\Re s > -C, \, \, |\Im s|<r$, which we will show is bounded above by $Cr^{1 + \delta}$, with $\delta$ the dimension of the Julia set.  We prove a weak, sublinear lower bound on the number of zeros in this region, an honest linear lower bound in logarithmic neighbourhoods of the imaginary axis, and conjecture that our upper bound is actually sharp.  We also obtain the upper bound $Cr^{1 + 2\alpha + \delta(1-\alpha)}$ for the number of zeros in more general regions $|\Re s | \leq |\Im s|^\alpha$ for $\alpha \in (0,1)$. \\
\indent Similar to \cite{Zworski1, GLZ}, we consider the dynamical system associated to a hyperbolic rational map $f$ when the Julia set $\JJ$ associated to this map is a Cantor-like totally disconnected set. \\
\indent We think of $\JJ$ as a subset of the sphere $\widehat{\cx} = \cx \bigcup \{\infty\}$ naturally identified with $\widehat{\reals}^2 = \reals^2 \bigcup \{ \infty\}$.  Then $|f'(z)|$ can be thought of as a map $\widehat{\reals}^2 \to \widehat{\reals}$, analytic in a neighbourhood of $\JJ$.  If $[f'(z)]$ is the holomorphic extension of $|f'(z)|$ to a map $\widehat{\cx}^2 = \cx^2 \bigcup \{\infty\} \to \widehat{\cx}$ define the transfer operator:
\begin{eqnarray}
\LL(s)u(z) = \sum\limits_{w \in f^{-1}(z)} \left[ f'(w)\right]^{-s} u(w).
\end{eqnarray}
\indent We will show in the course of this paper that $\LL(s)$ is a trace class operator on 
an appropriately chosen class of functions, and the Ruelle zeta function $Z$ will be 
defined as 
\begin{eqnarray} \label{zetafunction}
Z(s) = \det \left( I - \LL(s) \right).
\end{eqnarray}
We will prove the following bound of $Z$ in terms of $\delta$, the Hausdorff 
dimension of $\JJ$:
\begin{theorem}
\label{maintheorem}
Suppose $Z(s)$ is the zeta function defined by (\ref{zetafunction}) for the function $f$.  Then for any $C_0$, there exists $C_1$ such that for 
$\left|\Re s \right| \leq C_0$ we have
\begin{eqnarray}
\label{mainestimate}
\left|Z(s) \right| \leq C_1 \exp (C_1 \left| s \right|^{\delta}), \quad \delta = \dim \JJ
\end{eqnarray}
\noindent where $\delta$ is the dimension of the Julia set of $f$.
\end{theorem}
\indent This is the same result as in \cite{Zworski1}, but for more general Julia sets.  Using this and a dynamical formula for $Z(s)$ from Proposition \ref{dynamical-zeta}, we derive several estimates, both lower and upper bounds on the number of zeros in various regions.  Based on numerical evidence from \cite{Zworski1}, it appears the lower bounds are not optimal, and in closing we give an example to demonstrate the subtlety of this question and some of the problems in approaching it.

\section{Review of Julia sets}
\indent In this section we review a few classical results about the geometry of Julia 
sets that will be used later in the paper.  The interested reader should consult 
\cite{Falconer1, Falconer2} for further details. \\
\indent The Julia set $\JJ$ for a rational map can be defined to be the closure of the 
set of repelling periodic points, hence $\JJ$ is compact in the sphere.  It is easy to see \cite{Carleson} 
that $\JJ$ is backward and forward invariant: $\JJ = f(\JJ) = f^{-1}(\JJ)$, and in fact $f^p(\JJ) = \JJ$ for $p = 1,2, \ldots$.   We are interested in the case where $\JJ$ is disconnected (and hence totally disconneted).  The assumption that $f$ be hyperbolic means there exists an $n\geq 1$ such that $\inf \{|(f^n)'(z)|: z \in \JJ \} >1$.  In other words, some iterate of $f$ is expanding on the whole set.  A sometimes useful fact (see \cite{Urbanski}) is that a rational function is hyperbolic if and only if $\overline{\mathrm{PCV}(f)}\bigcap \JJ = \emptyset$, where $\mathrm{PCV}(f) = \bigcup_{n\geq 0} f^n(\mathrm{Crit}\,f)$ is the forward propagation of the set of critical points of $f$.  Note since $f$ is hyperbolic, we can replace $f$ with an appropriate iterate and assume that $f$ is strictly expanding near $\JJ$, so we will do this throughout. \\
\indent The most important properties of $\JJ$ are the those making it a ``cookie-cutter" set in the 
sense of \cite{Falconer2}.  Roughly speaking, this is to mean that a small neighbourhood intersected 
with $\JJ$ looks more or less like $\JJ$.  This is made precise in the following proposition:
\begin{proposition}
\label{cookie-cutter-prop}
$\JJ$ is a cookie-cutter set, that is, there exist constants $c >0$, 
$r_0>0$ such that for each $r<r_0$ and $z_0 \in \JJ$ there is a map $g: B(z_0,r) \to \widehat{\cx}$ such that 
$g(B(z_0,r) \bigcap \JJ) \subset \JJ$ satisfying 
\begin{eqnarray}
\label{cookie-cutter}
c^{-1}r^{-1}\left|z-w\right| \leq \left|g(z)-g(w)\right| \leq cr^{-1}\left|z-w\right|.
\end{eqnarray}
\end{proposition}
\indent To prove this, we will need the Koebe distortion theorem (see \cite{Carleson} for a proof):
\begin{lemma}[Koebe distortion theorem]
If $g$ is univalent (analytic and one-to-one) on the unit disk in $\cx$ with $g(0)=0$ and $g'(0)=1$, then 
\begin{eqnarray}
\label{Koebe}
\frac{1-\left|z\right|}{(1+\left|z\right|)^3} \leq \left|g'(z)\right| \leq 
\frac{1+\left|z\right|}{(1-\left|z\right|)^3}
\end{eqnarray}
\noindent and 
\begin{eqnarray}
\label{koebe}
\frac{\left|z\right|}{(1+\left|z\right|)^2} \leq \left|g(z)\right|\leq 
\frac{\left|z\right|}{(1-\left|z\right|)^2}
\end{eqnarray}
\end{lemma}
\indent We can get Proposition~\ref{cookie-cutter-prop} from this by a simple argument.  Since 
$\mathrm{Crit}\,f \bigcap \JJ = \emptyset$, there is $R>0$ so that for each $z_0 \in \JJ$, $f(z)$ is univalent on $B(z_0,R)$.  
This implies $f^n$ is also univalent, since $(f^n)'(z)= f'(f^{n-1}(z))\cdot f'(f^{n-2}(z)) \cdots 
f'(z)$.  We want to modify the estimate \ref{Koebe} to apply to a function $G$ univalent on a disk of 
radius $\delta>0$, say, with $G'(0)=M \neq 0$.  For $\zeta \in \left\{\left|\zeta\right|<1\right\}$, 
define $g(\zeta):= G(\delta \zeta)/M$.  Then $g$ satisfies the hypotheses of the lemma, and 
we now have:
\begin{eqnarray}
\left|M\right| \frac{1-\frac{1}{\delta}|z|}{\left(1 + \frac{1}{\delta}|z|\right)^3} \leq \left|G'(z)\right|
\leq \left|M\right| \frac{1+\frac{1}{\delta}|z|}{\left(1-\frac{1}{\delta}|z|\right)^4}
\end{eqnarray}
\noindent which, for $|z|<\delta/2$ yields
\begin{eqnarray}
\frac{|M|}{c} \leq \left| G'(z) \right|\leq |M|c
\end{eqnarray}
\noindent for some constant $c$.  The argument is finished by setting $r_0=R/2$ and noting that taking 
an appropriate iterate $G = f^n$ maps a ball $B(z_0,r)$ of radius $r<r_0$ centered at $z_0 \in \JJ$ into a larger fixed ball $B(0, S)$, say, with the property that $G(B(z_0,r) \bigcap \JJ) = \JJ$.  Thus there is a point $z_1 \in B(z_0,r)$ such that $S/(2r)\leq \left| G'(z_1) \right| \leq 2S/r$.  By conjugating with an appropriate M\"{o}bius transformation, we can assume $z_0 = z_1$, so that $C^{-1}r^{-1} \leq \left| G'(z) \right| \leq Cr^{-1}$ as claimed.

\section{The transfer operator on $L^2$ spaces}
It is more convenient to define the Ruelle transfer operator in terms of the inverse 
branches to $f$.  Suppose $f$ is an $m$ to $1$ function, and let $g_i(z)$ for $i = 1,2, \ldots,m$ be the branches of $f^{-1}$.  Now we interpret $\JJ$ as a subset of $\widehat{\reals}^2$ instead of 
$\widehat{\cx}$ and view $g_i:\widehat{\reals}^2 \to \widehat{\reals}^2$ real analytic and 
$\left|g_i' \right|:\widehat{\reals}^2 \to \widehat{\reals}$ analytic in a neighbourhood about $\JJ$.  
Then it is clear that both $g_i$ and $\left|g_i'\right|$ extend holomorphically to functions 
$g_i:\widehat{\cx}^2 \to \widehat{\cx}^2$ and $\left[ g_i' \right]:\widehat{\cx}^2 \to \widehat{\cx}$.  The Ruelle transfer operator 
can then be defined as 
\begin{eqnarray}
\label{Ldefine}
\LL(s)u(z) = \sum\limits_{i = 1}^m \left[g_i'(z)\right]^su(g_i(z)).
\end{eqnarray}
\indent We will show that with an appropriately chosen neighbourhood about $\JJ \subset \widehat{\cx}^2$ and an appropriately chosen class of functions $u$, $\LL$ is trace class.  We begin, as in 
\cite{Zworski1}, with a review of characteristic values of a compact operator $A:H_1 \to H_2$, 
where $H_j$'s are Hilbert spaces.  Define 
\begin{eqnarray*}
\| A \|=\nu_0(A)\geq \nu_1(A) \geq \cdots \geq \nu_l(A) \to 0
\end{eqnarray*}
\noindent to be the eigenvalues of $(A^*A)^{\frac{1}{2}} :H_1 \to H_1$.  The min-max principle 
shows that 
\begin{eqnarray}
\label{minmax}
\nu_l(A) = \min\limits_{{V \subset H_1} \atop {\mathrm{codim }V = l}} 
\max\limits_{{v \in V} \atop {\|v\|_{H_1}=1}} \|Av\|_{H_2}
\end{eqnarray}
\indent Now suppose $\left\{x_j \right\}_{j=0}^{\infty}$ is an orthonormal basis of $H_1$, 
then 
\begin{eqnarray}
\label{charestimate}
\nu_l(A) \leq \sum\limits_{j=l}^{\infty} \|Ax_j\|_{H_2}.
\end{eqnarray}
\noindent To see this we use $V_l = \mathrm{span }\left\{x_j\right\}_{j=l}^{\infty}$ in 
(\ref{minmax}): for $v \in V_l$ we have, by the Cauchy-Schwartz inequality, and the inequality 
$\|\cdot\|_{\ell_2} \leq \|\cdot\|_{\ell_1}$,
\begin{eqnarray*}
\|Av\|_{H_2}^2 = \left\| \sum\limits_{j=l}^{\infty}\langle v,x_j \rangle_{H_1}Ax_j \right\|^2 
\leq \|v\|_{H_1}^2 \left( \sum\limits_{j=l}^{\infty}\|Ax_j\|_{H_2}\right)^2,
\end{eqnarray*}
from which (\ref{minmax}) gives (\ref{charestimate}). \\
\indent We will also need the {\it Weyl inequality} (see \cite{traceideals}), which states 
that if $H_1 = H_2$ and $\lambda_j(A)$ are the eigenvalues of $A$,
\begin{eqnarray}
\left| \lambda_0(A) \right| \geq \left| \lambda_1(A) \right|\geq \cdots \geq \left| \lambda_l(A) 
\right| \to 0,
\end{eqnarray}
\noindent then for any $N$,
\begin{eqnarray}
\prod\limits_{l=0}^N(1 + \left|\lambda_l(A)\right|) \leq \prod\limits_{l=0}^N(1+\left|\nu_l(A)\right|).
\end{eqnarray}
\indent If $A$ is trace class, i.e. if $\sum_l\nu_l(A) < \infty$, then the determinant 
\begin{eqnarray*}
\det(I+A):=\prod\limits_{l=0}^{\infty}(1+\lambda_l(A)),
\end{eqnarray*}
\noindent is well defined and
\begin{eqnarray}
\label{nu-det-estimate}
\left|\det(I+A)\right| \leq \prod\limits_{l=0}^{\infty}(1+\nu_l(A)).
\end{eqnarray}
\indent We also need the following standard inequality about characteristic values (see 
\cite{traceideals})
\begin{eqnarray}
\label{nu-estimate}
\nu_{l_1+l_2+1}(A+B) \leq \nu_{l_1+1}(A) + \nu_{l_2+1}(B)
\end{eqnarray}
\noindent Lastly we finish with an obvious equality: suppose $A_j:H_{1j}\to H_{2j}$ and 
we form $\bigoplus_{j=1}^JA_j:\bigoplus_{j=1}^JH_{1j} \to \bigoplus_{j=1}^J H_{2j}$, then 
\begin{eqnarray}
\label{nu-direct-sum}
\sum\limits_{l=0}^{\infty}\nu_l\left( \bigoplus\limits_{j=1}^JA_j\right) = 
\sum\limits_{j=1}^J \sum\limits_{l=0}^{\infty}\nu_l(A_j).
\end{eqnarray}
\indent Now we want to define the Hilbert space we will be working with.  For 
$D \subset \widehat{\cx}^2$ open, let $H^2(D):= \left\{ u \mbox{ holomorphic in } D : \int_D 
|u(z)|^2dm(z) < \infty \right\}$.  We can take for $D$ disjoint neighbourhoods of $\JJ_i = g_i(\JJ)$ for 
$i=1,2, \ldots,m$ and we get the following:
\begin{proposition}
\label{trace-class-prop}
Suppose that $\LL(s): H^2(D) \to H^2(D)$ is defined by (\ref{Ldefine}), with $g_i$ the $m$ inverse branches of $f$ hyperbolic rational.  
Then for all $s \in \cx$, $\LL(s)$ is trace class and 
\begin{eqnarray}
\left| \det(I-\LL(s))\right| \leq C\exp(C|s|^3)
\end{eqnarray}
\noindent for some constant $C$.
\end{proposition}
\begin{proof}
We write $H^2(D) = \bigoplus\limits_{j=1}^{m}H^2(D_j)$ and $\LL(s) = \bigoplus\limits_{i,j = 1}^{m}\LL_{ij}(s)$, where 
\begin{eqnarray}
\label{lij}
\LL_{ij}(s)u(z):=\left[g_i'(z)\right]^su(g_i(z)), \,\,\, z\in D_j.
\end{eqnarray}
\indent Note that from (\ref{nu-estimate}) and (\ref{nu-direct-sum}) we certainly have 
\begin{eqnarray*}
\nu_k(\LL(s)) \leq m^2\max\limits_{1\leq i,j \leq m}\nu_{\left[ k/2m \right]}(\LL_{ij}(s)).
\end{eqnarray*}
Let $r_0>0$ be the minimum radius for which $|Dg_i(z)|<1$, $i=1,2, \ldots, m$, on a ball of radius $r_0$ centered at 
a point of $\JJ$.  Let 
\begin{eqnarray*}
U = \bigcup\limits_{i=1}^m U_i := \bigcup\limits_{i=1}^m\left\{\JJ_i+ B_{\widehat{\cx}^2}(0, r)\right\} 
\end{eqnarray*}
 for $r < r_0/2$.  Let $M = \max_{\overline{U}}|Dg_i(z)| 
<1$, and pick for $D$ a finite cover of $\JJ$, $D = \bigcup_{i=1}^m D_i$ made up of balls of 
radius $r$ centered at points of $\JJ$ as above, so that for each $z \in \JJ$, $d_{\widehat{\cx}^2}(z, \partial D) \geq \frac{1+M}{2}r$, and $D_i$ covers $\JJ_i = g_i(\JJ)$.  Then for any point $z \in D_i$, $|z-w|<r$ for some $w \in \JJ$, so $|g_j(z)-g_j(w)| 
\leq Mr$, so that $d_{\widehat{\cx}^2}(g_j(D_i), \partial D_j) \geq \left(\frac{1+M}{2}-M\right)r >0$.  
Lemmas \ref{char-lemma1} and \ref{char-lemma2} together with the estimate $\left| [g_i'(z)]^s\right| \leq e^{C|s|}$ 
now give for some $C_1$:
\begin{eqnarray*}
\nu_l(\LL_{ij}(s)) \leq C_1e^{C|s| -l^{1/2}/C_1}
\end{eqnarray*}
\indent With this in hand, we see that (\ref{nu-det-estimate}) implies 
\begin{eqnarray*}
\det(I-\LL(s)) \leq \prod\limits_{l=0}^{\infty}\left(1 + Ce^{C|s|-l^{1/2}/C}\right) \leq C e^{C^5|s|^3}
\end{eqnarray*}
\noindent so that $\LL(s)$ is trace class as claimed.  To finish with the proposition, we need the following 
two lemmas, taken almost directly from \cite{GLZ}.
\end{proof}
\begin{lemma}
\label{char-lemma1}
Let $\rho \in (0,1)$ and $R^{\rho}: H^2(B_{\cx^2}(0,1)) \to H^2(B_{\cx^2}(0,\rho))$ induced by the restriction 
map of $B_{\cx^2}(0,1)$ to $B_{\cx^2}(0,\rho)$.  Then for any $\widetilde{\rho} \in (\rho, 1)$ there exists a constant 
$C$ such that
\begin{eqnarray*}
\nu_l(R^{\rho}) \leq C \widetilde{\rho}^{l^{1/2}}
\end{eqnarray*}
\end{lemma}
\begin{proof}
Using (\ref{charestimate}) with the standard basis $(x_{\alpha})_{\alpha \in \mathbb{N}^2}$ for 
$H^2(B_{\cx^2}(0,1))$ given by
\begin{eqnarray}
x_{\alpha}(z) = c_{\alpha}z_1^{\alpha_1}z_2^{\alpha_2}, \quad \int_{B_{\cx^2}(0,1)} |x_{\alpha}(z)|^2 dz=1, \,\,\,
\alpha \in \mathbb{N}^2
\end{eqnarray}
\noindent for which we have 
\begin{eqnarray*}
\left\|R^{\rho}(x_{\alpha})\right\|^2 = \int_{B_{\cx^2}(0,\rho)} |x_{\alpha}(z)|^2 dz=\rho^{2|\alpha|+4}.
\end{eqnarray*}
\noindent The number of $\alpha$'s for which $|\alpha| \leq m$ is bounded by $(m+1)^2$, so by (\ref{charestimate}) 
\begin{eqnarray*}
\nu_l(R^{\rho}) \leq \sum\limits_{|\alpha| \geq l}\rho^{|\alpha|+2} \leq C \sum\limits_{k \geq l^{1/2}}(k+1)^2 \rho^k 
\leq C \widetilde{\rho}^{l^{1/2}}.
\end{eqnarray*}
\end{proof}
\begin{lemma}
\label{char-lemma2}
Suppose $\Omega_j \subset \cx^2$, $j=1,2$ are open sets and $\Omega_1=\bigcup_{k=1}^K B_{\cx^2}(z_k,r_k)$.  
Let $g$ be a holomorphic mapping defined on a neighbourhood $\widetilde{\Omega}_1$ of $\Omega_1$ taking values 
in $\Omega_2$ satisfying 
\begin{eqnarray*}
d_{\cx^2}(g(\Omega_1), \partial \Omega_2) > \frac{1}{C_0} >0, \quad 0 < \|Dg(z)\| < 1, \, \, z \in \widetilde{\Omega}_1.
\end{eqnarray*}
\noindent If 
\begin{eqnarray*}
A: H^2(\Omega_2) \to H^2(\Omega_1), \quad Au(z):= u(g(z)), \,\, z \in \Omega_1
\end{eqnarray*}
\noindent then for some $C_1$ depending only on the $r_k$'s, $K$, $d_{\cx^2}(g(\Omega_1), \partial \Omega_2)$, 
$\sup_{\widetilde{\Omega}_1}\|Dg(z)\|$, we have 
\begin{eqnarray*}
\nu_l(A) \leq C_1e^{-l^{1/2}/C_1}
\end{eqnarray*}
\noindent where $\nu_l(A)$'s are the characteristic values of $A$.
\end{lemma}
\begin{proof}
Define a new Hilbert space 
\begin{eqnarray*}
\mathcal{H}:= \bigoplus\limits_{k=1}^{K}H^2(B_k), \quad B_k = B_{\cx^2}(z_k, r_k),
\end{eqnarray*}
\noindent and a natural operator 
\begin{eqnarray*}
J: H^2(\Omega_1) \to \mathcal{H}, \quad (Ju)_k = u|_{B_k}
\end{eqnarray*}
\noindent We claim $J^*J: H^2(\Omega_1) \to H^2(\Omega_1)$ is invertible, with constants depending only on $K$.  To see this, note that for any $u \in H^2(\Omega_1)$,
\begin{eqnarray*}
\left\|u \right\|^2 \leq  \langle Ju,Ju \rangle_{\mathcal{H}} \leq K \left\|u\right\|^2.
\end{eqnarray*}
\noindent Hence 
\begin{eqnarray}
\label{JJ-estimate2}
\left\| J^*Ju \right\| = \langle Ju, JJ^*Ju \rangle_{\mathcal{H}} \leq K \left\|J^*Ju \right\|^2 \leq K^2 \left\| u \right\|^2
\end{eqnarray}
\noindent and
\begin{eqnarray}
\label{JJ-estimate1}
\left\| u \right\|^2 \leq \langle J^*J u, u \rangle_{H^2(\Omega_1)} \leq \left\| J^*Ju \right\| \left\| u \right\| 
\end{eqnarray}
\noindent for any $u \in H^2(\Omega_1)$.   The estimate (\ref{JJ-estimate2}) implies $J^*J$ is bounded, while the estimate (\ref{JJ-estimate1}) implies $J^*J$ is one-to-one.  Since any one-to-one self-adjoint operator is also onto, $J^*J$ is invertible, and furthermore,
\begin{eqnarray*}
\frac{1}{K^2} \left\|u\right\|^2 \leq \left\|(J^*J)^{-1}u\right\| \leq \|u\|^2.
\end{eqnarray*}
\indent Thus we calculate,
\begin{eqnarray*}
\nu_l(A) = \nu_l((J^*J)^{-1}J^*JA) \leq \left\|(J^*J)^{-1}\right\| \left\|J^*\right\| \nu_l(JA).
\end{eqnarray*}
\noindent Note then that 
\begin{eqnarray*}
\nu_{l}(JA) \leq K \max\limits_{1 \leq k \leq K} \nu_{[l/K]}(A_k), 
\end{eqnarray*}
\noindent where 
\begin{eqnarray*}
A_k: H^2(\Omega_2) \to H^2(B_k), \,\,\, A_ku(z) = u(g_k(z)), \,\,\, g_k= g|_{B_k}
\end{eqnarray*}
\noindent In order to estimate the characteristic values for $A_k$, note we can extend $g_k$ to a larger 
ball in $\widetilde{\Omega}_1$, $\widetilde{B}_k$ such that the image of its closure is still in $\Omega_2$.  
That gives us the operators $R_k: H^2(\widetilde{B}_k) \to H^2(B_k)$, $R_ku = u|_{B_k}$, and $\widetilde{A}_k$ 
defined similarly to $A_k$ with $B_k$'s replaced with $\widetilde{B}_k$'s.  Now we have $A_k = R_k\widetilde{A}_k$ 
which implies 
\begin{eqnarray*}
\nu_l(A_k) \leq \left\|\widetilde{A}_k\right\|\nu_l(R_k).
\end{eqnarray*}
\noindent Lemma \ref{char-lemma1} gives $\nu_l(R_k) \leq C_2 e^{-l^{1/2}/C}$ which completes the proof.
\end{proof}

\section{Estimates in terms of the dimension of $\JJ$}
In order to prove Theorem \ref{maintheorem} we need a few more important facts.  Recall that the diameter 
of a set $E$ is defined as $\diam(E) = \sup\left\{|x-y|:x,y\in E \right\}$.
\begin{proposition}
\label{diameter-proposition}
Let $\JJ \in \widehat{\cx}$ be the Julia set for $f$ hyperbolic rational.  Then there exist constants $K = K(c)$ and 
$\delta_0$ such that for $\delta < \delta_0$ the connected components of $\JJ + \overline{B}_{\widehat{\cx}}(0, \delta)$ 
have diameter at most $K \delta$.
\end{proposition}
\begin{proof}
Let $c$ and $r_0$ be as in Proposition \ref{cookie-cutter-prop}.  Since $\JJ$ is totally disconnected, there exists $\epsilon_0 > 0$ such that $\widehat{\JJ} = \JJ + B(0, \epsilon_0)$ has more than one connected component, and every connected component of $\widehat{\JJ}$ has diameter at most $(4c)^{-1}$.  Then we apply Proposition \ref{cookie-cutter-prop} with $r = c \delta \epsilon_0^{-1}$, with $\delta \leq \delta_0 < r_0 \epsilon_0 c^{-1}$.  The function $g$ guaranteed by Proposition \ref{cookie-cutter-prop} takes points in $\JJ$ to points in $\JJ$, so if $z \in \JJ$, $g(B(z, \delta)) \subset B(g(z), \epsilon_0)$.  Thus a connected component of $\JJ + B(0, \delta)$ is mapped into a connected component of $\widehat{\JJ}$.  Now suppose $d(z,w) > r/2$.  Then
\begin{eqnarray*}
d(g(z)-g(w)) \geq \frac{1}{cr}d(z,w) \geq \frac{1}{2c} > \frac{1}{4c}
\end{eqnarray*}
\noindent so that $g(z)$ is in a different connected component from $g(w)$.  Hence $z$ and $w$ must have been in separate connected components, and we conclude the diameter of the connected component containing $z$ is at most $K \delta=r$.
\end{proof}
\indent We have a bound on the diameter of the connected components of $\JJ + \overline{B}_{\widehat{\cx}}(0, \delta)$, but 
eventually we will need to cover $\JJ$ by balls, uniformly finite in $\delta$ so that we may again apply Lemma 
\ref{char-lemma2}.  
\begin{lemma}
\label{finite-cover-lemma}
Suppose $D\subset \widehat{\cx}$ is a compact set with the property that all connected components of $E = D + 
B_{\widehat{\cx}}(0, \delta)$ have diameter bounded by $K\delta$.  Then for any $\lambda \in (0,1)$ and any connected component $E_i$ of $E$, there 
exists a cover $U_i = U_i(\delta) \subset E_i$ of $D_i = E_i \bigcap D$ by at most $K' = K'(\lambda)$ balls of radius $\delta$ centered 
at points of $D_i$ such that $d_{\widehat{\cx}} (z, \partial U_i) \geq \lambda \delta$ for $z \in D_i$.
\end{lemma}
\begin{proof}
Let $l = (1-\lambda)/ \sqrt{2}$.  If $E_i$ is a connected component of $E$, then it fits in a closed ball of diameter $K\delta$ by hypothesis.  A ball of diameter $K\delta$ is contained in a closed cube $Q$ of sidelength $K\delta$, which can be covered by $K'(\lambda)$ closed cubes of side length $l\delta$ by simply starting at one corner of $Q$ and covering it with cubes $\{q_k\}$ of sidelength $l\delta$ intersecting only on their boundaries.  For each $k$, if $D \bigcap q_k \neq \emptyset$, select any point $p_k \in D \bigcap q_k$; if the intersection is empty, select nothing.  Then set $U_i = \bigcup\limits_{k=1}^{K'} B(p_k,\delta)$.  A simple calculation gives for any $z \in D_i$, $z \in q_k$ for some $k$, giving
\begin{eqnarray*}
d_{\widehat{\cx}}(z, \partial U_i) \geq d_{\widehat{\cx}}(p_k, \partial U_i) - d_{\widehat{\cx}}(p_k, z) \geq \lambda \delta
\end{eqnarray*}
\noindent since $z \in q_k$ implies, in particular, that $d_{\widehat{\cx}}(p_k, z) \leq \sqrt{2}l\delta$.
\end{proof}
\begin{remark}
It is clear that Lemma \ref{finite-cover-lemma} extends to $\widehat{\cx}^n$ with constants depending only on $K$ and the dimension.
\end{remark}
\begin{proof}[Proof of Theorem \ref{maintheorem}]
As in \cite{Zworski1} choose $h =|s|^{-1}$ where $|s|$ is large, but $|\Re (s)|$ is bounded.  
Now viewing $\JJ$ as a subset of $\widehat{\reals}^2$ instead of $\widehat{\cx}$, form $\widetilde{\JJ}(h) = \JJ + B_{\widehat{\cx}^2}(0, h)$.  Proposition \ref{diameter-proposition} tells us the diameter of each 
connected component of $\widetilde{\JJ}(h)$ has diameter less than $K h$.  Since $g_i$ (now thought of as the holomorphic extension $g_i : \widehat{\cx}^2 \to \widehat{\cx}^2$) is a contraction near $\JJ$, there 
is some $h_0$ so that if $h < h_0$, $M = \max |Dg_i(z)|<1$ for $z$ in the closure of $\JJ + B_{\widehat{\cx}^2}(0, h)$.  Let 
$\beta = \frac{1}{2} (M + 1)$, and suppose there are $P(h)$ connected components of $\widetilde{\JJ}(h)$.  Using Lemma 
\ref{finite-cover-lemma} we can pick a subcover $U(h)= \bigcup_{i=1}^m U_i(h)$ of at most $K' P(h)$ balls contained in 
$\widetilde{\JJ}$ and centered at points of $\JJ$ satisfying 
\begin{eqnarray*}
d_{\widehat{\cx}^2}(z, \partial U) \geq \beta h, \quad z \in \JJ.
\end{eqnarray*}
\noindent Since any point of $U$ is within $h$ of some $z \in \JJ$ and we know $g_i: \JJ_j \to \JJ_i$, 
\begin{eqnarray*}
d_{\widehat{\cx}^2}(g_i(U_j(h)), \partial U_i(h)) \geq (\beta -M)h > C^{-1}h
\end{eqnarray*}
\noindent for some constant $C$ independent of $h$.  \\
\indent It is classical that the Hausdorff measure of the Julia set is finite (see \cite{Urbanski} and the 
references therein) so that $P(h) = \mathcal{O}(h^{-\delta})$.  We write $\LL(s)$ as a sum of $m^2$ operators 
$\LL_{ij}(s)$ as before, $\LL_{ij}(s) : H^2(U_i) \to H^2(U_j)$, but now we have a better bound on the weight 
independent of $h$.  Since $[g_i'(z)]: \widehat{\reals}^2 \to \widehat{\reals}$,
\begin{eqnarray}
\left| [g_i'(z)]^s \right| & \leq & C \exp(|s||\arg \, [g_i'(z)]| ) \nonumber \\
&\leq & \exp\left(C_1 |s| \left(|\Im (z_1)|^2 + |\Im (z_2)|^2\right)^{1/2}\right) \label{g'-estimate} \\
& \leq & C_2,  \quad z \in U_j(h) \nonumber.
\end{eqnarray}
\noindent Each $\LL_{ij}(s)$ is a sum of no more than $P(h)$ operators, each of which satisfies $\nu_l \leq C \alpha^{l^{1/2}/C}$ for some $0 < \alpha < 1$ by Lemma \ref{char-lemma2}.  Thus using again (\ref{nu-estimate}) and 
(\ref{nu-direct-sum}) we get the estimate
\begin{eqnarray}
\log \left| \det \left( I -\LL(s) \right) \right| \leq C P(h) = \mathcal{O}(h^{-\delta})
\end{eqnarray}
\noindent which is (\ref{mainestimate}).
\end{proof}

\section{Counting Zeros in Strips}
In this section we prove the following corollary to Theorem \ref{maintheorem}.  The methods used here are similar to those used in \cite{Zworski1} and \cite{JP}.
\begin{corollary}
\label{counting-zeros}
Let $m(s)$ denote the multiplicity of a zero of $Z(s)$ at $s$.  Then 
\begin{eqnarray}
\label{counting-zeros-estimate}
\sum \left\{ m(s) : r \leq \left| \Im s\right| \leq r+1, \,\,\,\, \Re s > -C_0 \right\} \leq C_1 r^{\delta}
\end{eqnarray}
\noindent where $\delta = \dim \JJ$.
\end{corollary}
In order to prove this corollary, we will need to bound $Z(s)$ away from zero for $\Re s \geq C_0$.  We do this by employing a dynamical formula for $Z(s)$ which is interesting in its own right.  For the development of this dynamical formula, we take $D_i$ to be $\widehat{\cx}^2$-balls containing $\JJ \subset \widehat{\reals}^2$.  We again view $f$ as a map $f: \widehat{\reals}^2 \to \widehat{\reals}^2$ and then extend to a holomorphic function $\widehat{\cx}^2 \to \widehat{\cx}^2$ and write $f$ for this extension whenever unambiguous.
\begin{proposition}
\label{dynamical-zeta}
For $\Re(s) \gg 0$, 
\begin{eqnarray}
\label{dynamical-zeta-formula}
\det (I -\LL(s)) = \exp \left( - \sum\limits_{n=1}^{\infty}\frac{1}{n} \sum\limits_{f^n(z)=z} \frac{[(f^n)'(z)]^{-s}}{\left|\det(I - \left( d \, \left(f^n\right)(z)\right)^{-1})\right|}\right).
\end{eqnarray}
\end{proposition}
\begin{proof}
For $|\lambda|$ sufficiently small, $\log (I-\lambda \LL(s))$ is well defined and 
\begin{eqnarray*}
\det(I-\lambda \LL(s)) = \exp \left( -\sum\limits_{n=1}^{\infty} \frac{\lambda^n}{n} \tr(\LL(s)^n) \right).
\end{eqnarray*}
\noindent In order to evaluate the traces, we write
\begin{eqnarray*}
\tr \LL(s)^n = \sum\limits_{(i_1, \ldots, i_{n+1})} \tr\left(\LL_{i_1 i_2}(s) \circ \cdots \circ \LL_{i_ni_{n+1}}(s) \right),
\end{eqnarray*}
\noindent where $\LL_{ij}(s)$ is given by (\ref{lij}).  If the target space is different from the domain space, there are no eigenvalues, so that 
\begin{eqnarray*}
\sum\limits_{(i_1, \ldots, i_{n+1})} \tr\left(\LL_{i_1 i_2}(s) \circ \cdots \circ \LL_{i_ni_{n+1}}(s) \right) = \sum\limits_{(i_1, \ldots, i_{n})} \tr\left(\LL_{i_1 i_2}(s) \circ \cdots \circ \LL_{i_ni_{1}}(s) \right).
\end{eqnarray*}
\noindent We have
\begin{eqnarray*}
\LL_{i_1 i_2}(s) \circ \cdots \circ \LL_{i_ni_1}(s)u(z) = \left[(g_{i_1} \circ \cdots \circ g_{i_n})'(z) \right]^{s} u(g_{i_1} \circ \cdots \circ g_{i_n}(z)),
\end{eqnarray*}
\noindent and Lemma \ref{pullback-lemma} in the appendix shows that
\begin{eqnarray*}
\tr \left( (g_{i_1} \circ \cdots \circ g_{i_n})^*\right) = \frac{1}{\left|\det (I - d \, (g_{i_1} \circ \cdots \circ g_{i_n})(z))\right|}
\end{eqnarray*}
\noindent which completes the proof once we put $\lambda = 1$.
\end{proof}
\begin{proof}[Proof of Corollary \ref{counting-zeros}]
Using (\ref{dynamical-zeta-formula}) it is clear that for $\Re s \geq C_0$ we have 
\begin{eqnarray*}
|[(f^n)'(z)]^{-s}| \leq C C_1^{-n \Re (s)}
\end{eqnarray*}
\noindent with $C_1 > 1$ since $z$ is a periodic repeller and $[(f^n)'(z)]$ is real on $\JJ$.  Then the convergence of the double series in (\ref{dynamical-zeta-formula}) is immediate and gives for $\Re s \geq C_0$, $Z(s) \geq 1/2$.  With $Z$ zero free for $\Re s \geq C_0$, an application of the Jensen formula shows the left hand side of (\ref{counting-zeros-estimate}) is bounded by
\begin{eqnarray*}
\sum \left\{ m(s) : |s -ir -C_0| \leq 2(C_2+C_0) \right\} & \leq & C_3 \max\limits_{|s-ir-C_0|   \leq  4C_2} \log |Z(s)| \\
& \leq & C \max\limits_{{|\Re s |  \leq  C_3} \atop {|s| \leq 4C_2 + r}} \log |Z(s)| \\
& \leq & C_1 r^{\delta}.
\end{eqnarray*}
\end{proof}

\section{Polynomial Neighbourhoods}
Suppose $\{\mu_j \}$ are the zeros of $Z(s)$ counted with multiplicity.  Let $0 < \alpha < 1$ and consider the region $R_\alpha = \{ |\Re s| \leq | \Im s |^\alpha, |s| \geq 1 \}$.  Let $N_\alpha(r) = \#\{\mu_j \in R_\alpha : |\mu_j| \leq r \}$.  We expect $N_\alpha(r)$ to be somewhere in between the upper bound in strips, $N_0(r) = \# \{\mu_j: \Re \mu_j > C_0, \,\, | \Im \mu_j | \leq r \} \leq Cr^{1+ \delta}$, and the global bound $N_1(r) = \# \{ |\mu_j| \leq r \} \leq Cr^3$.  The following theorem asserts we get the expected interpolation.  The techniques used in the proof of theorem should extend easily to the case of convex co-compact Schottky groups \cite{GLZ}, giving the same upper bound as in \cite{Zworski3}.
\begin{theorem}
\label{main-theorem-2aa}
With $\{\mu_j\}$, $\alpha$, and $N_\alpha$ defined as above, 
\begin{eqnarray}
N_\alpha(r) \leq C_\alpha r^{1 + 2\alpha + \delta(1-\alpha)}
\end{eqnarray}
\end{theorem}
\begin{proof}
In order to begin, we start, as in the proof of Theorem \ref{maintheorem} by constructing a cover of $\JJ$ by open sets.  For $h = |s|^{-1}$, let $B_\alpha(h)$ be the open ball  
\begin{eqnarray*}
B_\alpha(h) := B_{\cx^2}(0, h^{1-\alpha}) , 
\end{eqnarray*}
and set $\widetilde{\JJ}(h) = \JJ + B_\alpha(h)$.  We can pick a finite subcover of $\widetilde{\JJ}(h)$, $U(h) = \bigcup_1^m U_i(h)$, as before, so that 
\begin{eqnarray*}
d(g_i(U_j(h), \partial U_i(h)) \geq C_\alpha^{-1}h^{1-\alpha}
\end{eqnarray*}
\noindent for some constant $C_\alpha$ independent of $h$.  Write $\LL$ as a sum of $m^2$ operators $\LL_{ij}(s): H^2(U_i) \to H^2(U_j)$ as before.  Since $[g_i'(z)]: \widehat{\reals}^2 \to \reals$, if we take $s \in R_\alpha' := \{ |\Re s | \leq 5 |\Im s|^\alpha \}$,
\begin{eqnarray}
|[g_i'(z)]^s| & \leq & C \exp (C|\Re s||\log (|[g_i'(z)]|) + C |\Im s| |\arg [g_i'(z)] |) \nonumber \\
& \leq & C \exp \left( C |s|^{\alpha} + C |\Im s| \left( (\Im z_1)^2 + (\Im z_2)^2 \right)^{\frac{1}{2}} \right) \\
& \leq & C \exp (C |s|^\alpha). \nonumber
\end{eqnarray}
\indent Now each $\LL_{ij}(s)$ is a sum of no more than $P(h)= \O(h^{-\delta(1-\alpha)})$ operators, each of which has characteristic values $\{\nu_l\}$ satisfying $\nu_l \leq Ce^{C|s|^\alpha - l^{1/2}/C}$.  Then for $s \in R_\alpha'$,
\begin{eqnarray*}
|Z(s)| & \leq & \prod\left( 1 + C e^{C|s|^{\alpha}-l^{1/2}/C} \right)^{|s|^{\delta ( 1-\alpha)}} \\
& \leq & C e^{C^5|s|^{3 \alpha + \delta(1-\alpha)}}.
\end{eqnarray*}
\indent Next we restrict attention to zeros in the upper left quadrant.  Observe that for any $r_1 < r_2$, we can cover $R_\alpha \bigcap \{ r_1 \leq |s| \leq r_2 \}$ by boxes of width $2r^\alpha$ and height $r^\alpha$ with right bottom corner at $s = ir$, for $r_1 < r < r_2$.  If $n_B(r)$ is the number of such boxes needed to cover $R_\alpha$, clearly $r^{-\alpha}/2 \leq dn_B(r) \leq 2r^{-\alpha}$.  Each box $[-2r^\alpha,0] \times i[r, r+r^\alpha]$ can be covered by two discs $D_1(r) \bigcup D_2(r)$ with 
\begin{eqnarray*}
D_1(r) & = & D(-\frac{3}{2}r^\alpha+i(r + \frac{1}{2}r^\alpha), \frac{1}{\sqrt{2}}r^\alpha), \\ 
D_2(r) & = & D(-\frac{1}{2}r^\alpha + i(r+\frac{1}{2}r^\alpha), \frac{1}{\sqrt{2}}r^\alpha), 
\end{eqnarray*}
\noindent both of which fit inside of $R_\alpha'$ (see Figure \ref{r-alpha-fig}).  The Jensen formula tells us $N_D(r) :=\# \{ \mu_j : \mu_j \in D_1(r) \bigcup D_2(r)\} \leq C r^{3 \alpha + \delta(1-\alpha)}$.  Then 
\begin{eqnarray*}
N_\alpha(r) \leq \int_1^r N_D(s) dn_B(s) \leq Cr^{1 + 2 \alpha + \delta(1- \alpha)}
\end{eqnarray*}
\noindent as claimed.
\begin{figure}
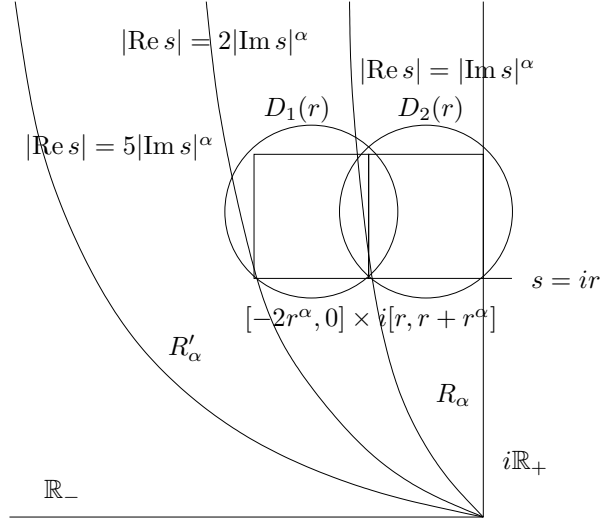

\caption{\label{r-alpha-fig} Regions used in the proof of Theorem \ref{main-theorem-2aa}}
\include{R-fig2}
\end{figure}
\end{proof}

\section{Lower Bounds on the Number of Zeros}
In order to prove lower bounds on the number of zeros, we will use extensively the dynamical formula (\ref{dynamical-zeta-formula}).  In light of \cite{CM} we see the series in Proposition \ref{dynamical-zeta} actually converges for all $\Re s > \delta$.  We will use this in the following proofs when we select our contours of integration.  Let $w(s) = Z(is + \delta)$, and suppose $\{\lambda_j\}_{j=1}^{\infty}$ are the zeros of $w$ counted with multiplicity.  Let $u_1(t) \in \mathcal{D}'(\reals_+)$ be the distribution 
\begin{eqnarray}
\label{u1}
u_1(t) : = \sum_j e^{it \lambda_j}
\end{eqnarray}
\noindent and let $u_2(t) \in \mathcal{D}'(\reals_+)$ be the distribution
\begin{eqnarray}
\label{u2}
u_2(t) := t e^{-\delta t} \sum_n \frac{1}{n} \sum_{f^n(z) = z}\frac{\delta_0(t-L_n(z))}{|\det(I-d\,(f^n)^{-1}(z))|}
\end{eqnarray}
\noindent where $\delta_0$ is the usual Dirac mass and 
\begin{eqnarray*}
L_n(z) = \log [(f^n)'(z)].
\end{eqnarray*}
\begin{lemma}
\label{distribution-id}
With $u_1$, $u_2$ as above, 
\begin{eqnarray}
\label{distribution-id-statement}
u_1(t) = u_2(t) 
\end{eqnarray}
in the sense of distributions on $\reals_+$.
\end{lemma}
\begin{remark} We use this distribution identity to make the presentation of the following proofs clear, and in order to quote directly Lemma \ref{sjoz-thm} below from \cite{SjoZ}.
\end{remark}
\begin{proof}
Let $w_\epsilon(s) = Z(is + \delta + \epsilon)$ for $\epsilon > 0$.  Then $w_\epsilon$ has a dynamical expansion for $\Im s < 0$, and if $\{\lambda_j^{\epsilon}\}$ are the zeros of $w_\epsilon$ counted with multiplicity, then $\Im \lambda_j^\epsilon > 0$ for each $j$.  In light of Proposition \ref{trace-class-prop}, we see $w(s)$ is an entire function of order $3$, hence the Weierstrass factorization gives
\begin{eqnarray}
\label{weierstrass}
\frac{d^3}{d\lambda^3} \left( \frac{w'(\lambda)}{w(\lambda)} \right) = -3! \sum_j \frac{1}{(\lambda_j -\lambda)^4}.
\end{eqnarray}
\noindent Now 
\begin{eqnarray*}
-3! \frac{1}{(\lambda_j - \lambda)^4} = i\frac{d^3}{d \lambda^3} \int_0^{\infty} e^{it(\lambda_j - \lambda)} dt = -\int_0^{\infty}t^3 e^{it(\lambda_j - \lambda)} dt
\end{eqnarray*}
\noindent provided $\Im \lambda_j >0$ and $\lambda$ is real.  Hence the right hand side of (\ref{weierstrass}) is 
\begin{eqnarray}
\label{u1-derivatives}
v_1^\epsilon(t) =  -\mathcal{F} \left( t^3\sum_j (e^{it \lambda_j^{\epsilon}})_+ \right)
\end{eqnarray}
\noindent and the right hand side of (\ref{weierstrass}) is 
\begin{eqnarray}
\label{u2-derivatives}
v_2^\epsilon (t) =  -\mathcal{F} \left(  t^4 e^{-(\delta + \epsilon) t} \sum_n  \frac{1}{n} \sum_{f^n(z)=z} \frac{\delta_0 ( t-L_n(z))}{|\det (I - d\,(f^n)^{-1}(z))|} \right),
\end{eqnarray}
\noindent where $\mathcal{F}$ denotes the usual Fourier transform.  Since both of these distributions are tempered and $\mathcal{F}$ is an isomorphism on $\mathcal{S}'$, we conclude 
\begin{eqnarray*}
u_1^\epsilon(t) = u_2^\epsilon(t)
\end{eqnarray*}
\noindent away from $0$, where 
\begin{eqnarray}
\label{u1-epsilon}
u_1^\epsilon(t):= \sum_j e^{it \lambda_j^{\epsilon}}
\end{eqnarray}
\noindent and 
\begin{eqnarray}
\label{u2-epsilon}
u_2^\epsilon(t) :=  t e^{-(\delta + \epsilon) t} \sum_n \frac{1}{n} \sum_{f^n(z)=z} \frac{\delta_0 ( t-L_n(z)}{|\det (I - d\,(f^n)^{-1}(z))|}.
\end{eqnarray}
Finally, integrating against a test function in $C_0^\infty(\reals_+)$ and sending $\epsilon$ to zero gives (\ref{distribution-id-statement}).
\end{proof}
\indent Next we use this distribution identity to count zeros in specific regions.

\subsection{Zeros in Strips}
We first recall some notation: we say $f(x)=\Omega (g(x))$ as $x \to \infty$ if there does not exist any constant $C$ for which $f(x) \leq C g(x)$ as $x \to \infty$.  That is, $f(x)$ cannot be controlled by $g(x)$ as $x \to \infty$.
\begin{theorem}
\label{main-theorem-2a}
Let $\delta$ be the dimension of the Julia set, $Z(s)$ the dynamical zeta function, and $\{\mu_j \}$ the zeros of $Z$ with multiplicity.  Then there exists $\epsilon_0 > 0$ such that for all $0 < \epsilon < \epsilon_0$, $\#\{ \mu_j : | \Im \mu_j | \leq r, \,\,\,  - C\epsilon^{-1} < \Re \mu_j  < \delta \} = \Omega (r^{1-\epsilon})$.
\end{theorem}
Observe Corollary \ref{counting-zeros} implies an upper bound on the number of zeros in strips:
\begin{eqnarray*}
\#\left\{ \mu_j : |\Im \mu_j| \leq r, -C < \Re \mu_j < \delta \right\} \leq Cr^{1+\delta}
\end{eqnarray*}
which suggests this lower bound is in fact not optimal.  Instead we have the following conjecture.
\begin{conjecture}
\label{lower-bound-conjecture}
There exists $\epsilon_0 > 0$ such that for all $0 < \epsilon < \epsilon_0$, there exists $0<C_\epsilon< \infty$ such that
\begin{eqnarray*}
\#\{ \mu_j : | \Im \mu_j | \leq r, \,\,\,  - \epsilon^{-1} < \Re \mu_j  < \delta \} \geq (C_{\epsilon}^{-1})r^{1+ \delta},
\end{eqnarray*}
\noindent where $\delta $ is the dimension of the Julia set.
\end{conjecture}
\begin{proof}[Proof of Theorem \ref{main-theorem-2a}]
Let $\P$ satisfy $\Ph \in C_0^{\infty}(\reals)$, $\Ph(0) = 1$, $\Ph \geq 0$ and $\supp \Ph \subset [-1,1]$, where $\Ph$ denotes the Fourier transform.  Define $\Phd(t) = \Ph(\gamma^{-1}(t-d)$ with $d \geq 1$, $\gamma \leq 1$, so that $\supp \Phd \subset [d -\gamma, d+\gamma] \subset \reals_+$.  Then 
\begin{eqnarray}
\label{ph-u1}
\langle u_1, \Phd \rangle = (2 \pi )^{-1/2}\sum_j \Pd (\lambda_j)
\end{eqnarray}
\noindent and 
\begin{eqnarray}
\label{ph-u2}
\langle u_2, \Phd \rangle = \sum_n \frac{1}{n}\sum_{f^n(z)=z}\frac{L_n(z) \exp(-L_n(z)\delta )}{|\det(I - d\, (f^n)^{-1}(z))|} \Phd(L_n(z))
\end{eqnarray}
\noindent with $L_n(z) = \log [(f^n)'(z)]$.  If $d$ is chosen near one of the $L_n(z)$s and $\gamma$ is small, (\ref{ph-u2}) is bounded from below by 
\begin{eqnarray}
\label{lower-bound-1}
C^{-1}d e^{-\delta d}.
\end{eqnarray}
\indent Next we deal with (\ref{ph-u1}).  Note $\Pd$ is an entire holomorphic function satisfying
\begin{eqnarray*}
\left| \Pd( \zeta ) \right| = \left| \gamma  \P ( \gamma \zeta ) \right| \leq C_M \gamma \exp\left( (\gamma - d) \Im \zeta \right) \left( 1 + |\gamma \zeta | \right)^{-M}
\end{eqnarray*}
\noindent for any $N$, $\Im \zeta \geq 0$ by the Paley-Wiener theorem \cite{hormander}.  \\
\indent Since $w$ is entire of order $3$, the Jensen formula gives $N_1(r) = \# \{ \lambda_j : |\lambda_j| \leq r \} \leq Cr^3$.  Then for $\kappa > 0$,
\begin{eqnarray}
\label{geq-kappa}
\lefteqn{\left| \sum_{\{ \Im \lambda_j \geq \kappa\}}  \Pd(\lambda_j) \right| } \nonumber \\
&\leq& C \gamma e^{(\gamma-d) \kappa} \int_0^\infty (1 + \gamma r)^{-M} r^2 dr \\
&\leq& C  \gamma^{-2} e^{(1-d)\kappa}.
\end{eqnarray}
Now assume 
\begin{eqnarray*}
N(\kappa,r) = \# \left\{ \lambda_j : \left| \Re \lambda_j \right| \leq r, \,\, \Im \lambda_j < \kappa \right\} \leq C_{\epsilon}(\kappa) r^{1-\epsilon}
\end{eqnarray*}
\noindent for some constant $C_\epsilon(\kappa)$.  Then
\begin{eqnarray} 
\label{leq-kappa}
\left| \sum_{\{ \Im \lambda_j < \kappa \} } \Pd(\lambda_j) \right| \leq C  \gamma \int_0^\infty (1 + \gamma r) ^{-M} N(\kappa, dr) \leq C\gamma^\epsilon.
\end{eqnarray}
\noindent Combining (\ref{geq-kappa}) and (\ref{leq-kappa}) we get that (\ref{ph-u2}) is bounded from above by 
\begin{eqnarray*}
C \gamma^{-2} e^{(1-d) \kappa} + C  \gamma^\epsilon.
\end{eqnarray*}
Hence we have the inequality
\begin{eqnarray*}
C^{-1}de^{-\delta d} \leq C  \gamma^{-2} e^{(1-d) \kappa} + C \gamma^\epsilon
\end{eqnarray*}
\noindent which yields a contradiction once we set $\gamma = e^{- \beta d}$ with $\beta > \delta / \epsilon $, and $C \epsilon^{-1} > \kappa > \delta +  2( \delta / \epsilon) $.
\end{proof}

\subsection{Lower Bounds in Logarithmic Neighbourhoods}
In this section, we use a theorem from \cite{SjoZ} to get improved lower bounds in logarithmic neighbourhoods of the real axis.  To this end, let $\Lambda = \{ \lambda_j \}$ be the set of zeros for $w(s) = Z(is+ \delta )$, and let $\Lambda_{\rho} = \{ \lambda_j : \Im \lambda_j < \rho \log |\lambda_j| \}$.  Let 
\begin{eqnarray*}
N_1(r) &=& \# \{\lambda_j : |\lambda_j| \leq r \}, \\  N_{\rho}(r) &=& \# \{ \lambda_j \in \Lambda_{\rho} : |\lambda_j| \leq r \}.  
\end{eqnarray*}
We know $N_1(r) = \O(r^3)$ from Proposition \ref{trace-class-prop}.  We use a slightly different test function for this development.  Let $\P \in C_0^{\infty}(\reals)$ satisfy $\supp(\P) \subset [-1,1]$, $\P(0) = 1$, $\Ph( \zeta ) \geq 0$ for all $\zeta \in \reals$ and for $d>1$, $\gamma <1$, set $\Pd = \P(\gamma^{-1}(t-d))$.  We will need the following lemma, taken directly from \cite{SjoZ}:
\begin{lemma}
\label{sjoz-thm}
Suppose $\{\lambda_j \} \subset \cx$ is a sequence of points such that $u(t) := \sum_j e^{it \lambda_j}$ belongs to $\mathcal{D}'(\reals_+)$.  Suppose for some $k \in \reals$ and fixed $d>0$ there is a constant $0<C< \infty$ such that 
\begin{eqnarray}
\label{sjoz-hyp}
\left| \widehat{( \P u )}(\lambda) \right| > \left( C^{-1} - o(1) \right) |\lambda|^k, \quad \quad \lambda \to \infty
\end{eqnarray}
\noindent for every $\P \in C_0^\infty(\reals_+)$ with sufficiently small support such that $\P(d) = 1$.  Then for every sufficiently small $\epsilon > 0 $ and $\rho > (n-k)/(d - \epsilon^2)$ we have: \\
\noindent a) If $k \geq 0$, then 
\begin{eqnarray}
\label{sjoz-bound}
N_\rho (r) > \left( \frac{1}{ C \pi ( k+1)} - o(1) \right)r^{k+1}
\end{eqnarray}
\noindent and moreover
\begin{eqnarray*}
\sum_{ \lambda_j \in \Lambda_\rho, \,\,\, |\Re \lambda_j| \leq r } e^{-(d-\epsilon)\Im \lambda_j} > \left( \frac{1}{ C \pi ( k+1)} - o_\epsilon(1) \right)r^{k+1}.
\end{eqnarray*}
\noindent b) If $k < 0$, then for every $\eta > 0$, there is $r(\eta)>0$ such that 
\begin{eqnarray*}
N_\rho(r) > r^{1 - \eta} \,\,\,\, \mathrm{if} \,\,\,\, r > r(\eta).
\end{eqnarray*} 
\end{lemma}
For a proof see \cite{SjoZ}.  We will use the first part of the lemma to deduce the following ``honest'' linear lower bound:
\begin{corollary}
\label{linear-lower-bound-1}
For $\Pd$ and $\rho$ as above, 
\begin{eqnarray}
\label{corollary-bound-statement}
N_\rho(r) \geq (C^{-1} - o_\gamma(1))r.
\end{eqnarray}
\end{corollary}
\begin{proof}
Using the distribution identity (\ref{distribution-id}), we see $u_1$ is of the correct form to apply Lemma \ref{sjoz-thm}.  It remains then only to verify that (\ref{sjoz-hyp}) holds with $k = 0$.  If $f(z)$ is hyperbolic rational, we can replace $f$ with an appropriate iterate so that $f'(z)>1$ on $\JJ$.  Then $ n \log(A) \leq L_n(z) \leq n \log(B)$ for all $n$, $f^n(z) = z$, where $A = \min_{\JJ} f'(z)$ and $B = \max_{\JJ} f'(z)$.  Since there are precisely $m^n$ discrete orbits for each $n$ and the $L_n(z) \to \infty$, if we fix $n$ we can find $\gamma_n$ small enough and $\ell$ close to $n\log(AB)^{1/2}$ so that $\ell = L_n(z)$ for at least one orbit and  
\begin{eqnarray*}
\varphi_{ \gamma, \ell }(L_n(z)) =
\left\{ \begin{array}{l} 1 \,\,\,  \mathrm{if} \,\,\, L_n(z) =  \ell \\
			 0 \,\,\, \mathrm{otherwise}.  \end{array} \right.
\end{eqnarray*}
Then for this $\Pd$, we calculate for $u_1, u_2$ defined above:
\begin{eqnarray}
\label{pd-u2-2}
|\widehat{(u_1 \Pd)}(\lambda) | & = & |\widehat{(u_2 \Pd)}(\lambda) | \nonumber \\
& = & \left| \sum_n \frac{1}{n} \sum_{f^n(z)=z} \left(\frac{L_n(z) e^{-\delta L_n(z)} \Pd(L_n(z))}{|\det(I - d\, (f^n)^{-1}(z))|} e^{-iL_n(z) \lambda} \right) \right| \nonumber \\
& = & \label{lower-bound-equality1}   \sum_n \frac{1}{n} \sum_{{f^n(z)=z} \atop {L_n(z) = \ell}} \frac{\ell e^{-\delta  \ell} }{|\det(I - d\, (f^n)^{-1}(z))|}    \\
& > & \label{lower-bound-inequality1} C^{-1} 
\end{eqnarray}
\noindent with (\ref{lower-bound-equality1}) and (\ref{lower-bound-inequality1}) holding because the sums are finite and all terms are positive.  Thus $u_1$ satisfies the hypotheses of Lemma \ref{sjoz-thm} with $k=0$ and (\ref{sjoz-bound}) gives (\ref{corollary-bound-statement}).
\end{proof}

\section{Final Comments}
Experimental evidence in \cite{Zworski1} suggests Conjecture \ref{lower-bound-conjecture} is true.  However, as is common with this type of estimates, sharp lower bounds have remained elusive.  In order to illustrate the subtlety of this question, we will look at the following example.  \\
\indent Assume for simplicity that $f(z) = z^2+c$ for $c$ real, $c<-2$, and that $A/B$ is irrational, with $1< A=\min_{\JJ} |f'(z)|<B = \max_{\JJ} |f'(z)|$ as before.  Then $\JJ$ is a Cantor-like set in the real line and all the proofs above go through by complexifying to $\cx$ instead of $\cx^2$.  In the proof of Corollary \ref{linear-lower-bound-1}, we stated that the distribution of the $L_n(z)$s is Gaussian with concentration at $\log (AB)^{1/2}$.  This suggests a simple model for the zeta function.  With $A$ and $B$ as above, we model the distribution of the $L_n(z)$s in the following fashion.  We write $L_n(z) = kl_1 + (n-k)l_2$ with multiplicity ${n \choose k}$, where we have set $l_1= \log A$ and $l_2 = \log B$.  Using (\ref{dynamical-zeta-formula}) as a basis, we calculate
\begin{eqnarray*}
\lefteqn{-\sum \frac{1}{n} \sum_{f^n(z)=z} \frac{\exp \left( -sL_n(z) \right)}{\left(1-\exp ( -L_n(z))\right)} =}  \\
& = & -\sum_n \frac{1}{n} \sum_{f^n(z) = z} \sum_k \exp \left( -(s+k)L_n(z) \right) \\
& = & -\sum_n \frac{1}{n} \sum_k \sum_m {n \choose m} \left( \exp (-(s+k)l_1)\right)^m \left( \exp (-(s+k)l_2 \right)^{n-m} \\ 
& = & -\sum_n \frac{1}{n} \sum_k \left( \exp (-(s+k)l_1) + \exp (-(s+k)l_2)\right)^n \\
& = & -\sum_k \log \left( 1 - e^{-(s+k)l_1}-e^{-(s+k)l_2} \right) \\
& = & \log \prod_k \left( 1 - e^{-(s+k)l_1}-e^{-(s+k)l_2} \right) 
\end{eqnarray*}
\noindent so we set  
\begin{eqnarray}
\widetilde{Z}(s) = \prod_k\left( 1 - A^{-(s+k)}-B^{-(s+k)} \right).
\end{eqnarray} 
This model shares some important features with $Z(s)$.  First, it has one zero at $s= \delta$, where $\delta$, solving $A^\delta + B^\delta = 1$, is the ``dimension''.  Second, it is easy to see that, since $A/B$ is irrational, there are no other zeros on $\Re s = \delta$.  However, if $\Re s > -C$, we can take \begin{eqnarray*}
\widetilde{Z}(s) \sim C\prod_{k=0}^K \left( 1 - A^{-(s+k)}-B^{-(s+k)} \right)
\end{eqnarray*}
\noindent for some $K$.  Then as $|s| \to \infty$, $\widetilde{Z}(s) \sim Ce^{C|s|}$, whence the number of zeros in $\{\Re s > -C, |\Im s| \leq r \}$ grows linearly.

\appendix
\section{}
The following lemma is widely known in the literature, but we include a proof of the general result here for completeness.
\begin{lemma}
\label{pullback-lemma}
Suppose $\Omega \subset \reals^2$ is an open, bounded, convex domain, and $f: \Omega \to \Omega$ is an analytic contraction obtained from a holomorphic function on $\cx$ identified with $\reals^2$.  Let $\tilde{f}: \widetilde{\Omega} \to \widetilde{\Omega}$ denote the extension of $f$ to a holomorphic contraction on a bounded, pseudoconvex domain $\widetilde{\Omega} \subset \cx^2$.  Suppose $z_1$ is the unique fixed point of $\tilde{f}$.  Then the pullback operator $\tilde{f}^* : H^2 ( \tilde{f}(\widetilde{\Omega})) \to H^2(\widetilde{\Omega})$ has trace 
\begin{eqnarray}
\label{pullback-formula-statement}
\tr \tilde{f}^* = \frac{1}{\left|\det \left( I - d \, \tilde{f}(z_1)\right)\right|}.
\end{eqnarray}
\end{lemma}
We first prove this result in the case $\widetilde{\Omega}$ is a ball.
\begin{lemma}
\label{pullback-lemma'}
Suppose $f: B_{\reals^2}(z_0,r) \to B_{\reals^2}(z_0, r')$ is an analytic contraction obtained from a holomorphic function on $\cx$ identified with $\reals^2$, and let $\tilde{f} :B_{\cx^2}(z_0,r) \to B_{\cx^2}(z_0, r')$ be the holomorphic extension of $f$ to $\cx^2$.  If $z_1$ is the unique fixed point of $\tilde{f}$, then the pullback by $\tilde{f}$, $\tilde{f}^*:H^2(B_{\cx^2}(z_0, r')) \to H^2(B_{\cx^2}(z_0, r))$ has trace
\begin{eqnarray*}
\tr \tilde{f}^* = \frac{1}{\left|\det \left( I - d \, \tilde{f}(z_1)\right)\right|}.
\end{eqnarray*}
\end{lemma}
\begin{proof}
Without loss of generality, $z_0 = 0$, and $\tilde{f}: B_{\cx^2}(0,1) \to B_{\cx^2}(0, \rho)$ for some $\rho \in (0,1)$.  Since the group $SU_{\cx}(2,1)$ acts transitively on the unit ball in $\cx^2$, by composing with appropriate M\"obius transformations we may also assume $z_1=0$ (see \cite{ahlfors}).  We first consider $f: B_{\reals^2}(0,1) \to B_{\reals^2}(0, \rho)$.  The assumption that $f$ is obtained from a holomorphic function on $\cx$ means for $z \in \cx$, 
\begin{eqnarray*}
d \, f(0) z = (a + ib)(x+iy)=(ax-by) + i(bx + ay)
\end{eqnarray*}
\noindent for some $a,b \in \reals$.  But this implies $d\, f_{\reals^2}(0)$ and hence $d\, \tilde{f}(0)$ has the very special form
\begin{eqnarray*}
d\, \tilde{f}(0) = \left(\begin{array}{rr} a & -b \\ b & a \end{array} \right), \quad \quad a,b \in \reals.
\end{eqnarray*}
\noindent Thus $d \, \tilde{f}(0)$ is always diagonalizable.  Note then that if 
\begin{eqnarray*}
d \, \tilde{f}(0) \left(\begin{array}{c}z_1 \\ z_2 \end{array} \right) = \left( \begin{array}{c} az_1-bz_2 \\ bz_1 + az_2 \end{array} \right),
\end{eqnarray*}
\noindent then the change of variables 
\begin{eqnarray}
\lefteqn{ \left( \begin{array}{c} w_1 \\ w_2 \end{array} \right) = A \left( \begin{array}{c} z_1 \\ z_2 \end{array} \right) = } \nonumber\\  & = & \frac{1}{a^2+b^2}\left( \begin{array}{cc} a(a-ib) & -b(a-ib) \\ b(a+ib) & a(a+ib) \end{array} \right) \left( \begin{array}{c} z_1 \\ z_2 \end{array} \right) \label{change-variables}
\end{eqnarray}
\noindent makes $d\, \tilde{f}(0)$ diagonal, and further, $\det A = 1$. \\
\indent We have an orthonormal basis for $H^2(B_{\cx^2}(0,1))$ in the form $\{ c_\alpha z^{\alpha} \}_{\alpha \in \mathbb{N}^2}$ for constants $c_\alpha$.  We can use the Bergman kernel to write the kernel for the pullback operator on $H^2(B_{\cx^2}(0,1))$, 
\begin{eqnarray*}
K_{\tilde{f}^*}(z, s) = \sum_{\alpha} |c_\alpha|^2 (\tilde{f}(z))^{\alpha} \overline{s}^{\alpha},
\end{eqnarray*}
\noindent so that for each $u \in H^2(B_{\cx^2}(0,1))$,
\begin{eqnarray*}
\tilde{f}^*u(z) = \int_{B_{\cx^2}(0,1)} K_{\tilde{f}^*}(z,s) u(s)dm(s).
\end{eqnarray*}
\noindent Here $dm(s)$ denotes the usual Lebesgue measure on $\cx^2$.  We will use the change of variables (\ref{change-variables}) and the fact that $\tilde{f}^*$ is trace class to exchange the integral and sum in the following to get: 
\begin{eqnarray*}
\tr \tilde{f}^* &=& \int_{B_{\cx^2}(0,1)} K_{\tilde{f}^*}(z,z) dm(z) \\
& = & \int_{B_{\cx^2}(0,1)} \sum_{\alpha}|c_{\alpha}|^2 (\tilde{f}(z))^{\alpha} \overline{z}^\alpha dm(z) \\
& = & \int_{B_{\cx^2}(0,1)} \sum_{\alpha} |c_{\alpha}|^2 (d \, \tilde{f}(0) z + \O(|z|^2))^\alpha \overline{z}^\alpha dm(z) \\
& = & \int_{B_{\cx^2}(0,1)} \sum_{\alpha} |c_{\alpha}|^2 \left( \left( \begin{array}{cc} a+ib & 0 \\ 0 & a-ib \end{array} \right) w + \O(|w|^2) \right)^{\alpha} \overline{w}^\alpha dm(w) \\
& = & \sum_{\alpha} (a+ib)^{\alpha_1} (a-ib)^{\alpha_2} \\
& = & \left| \det (I - d\, \tilde{f}(0)) \right|^{-1}.
\end{eqnarray*}
\end{proof}
\begin{proof}[Proof of Lemma \ref{pullback-lemma}]
Let $B$ be the largest open ball with center at $z_1$ so that $B \subset \widetilde{\Omega}$.  Since $\tilde{f}$ is a contraction and we can always replace $f$ with an appropriate iterate if necessary, we may assume without loss of generality that $\tilde{f}(\widetilde{\Omega}) \subset B$.  Now suppose $u$ is a generalized eigenfunction of $\tilde{f}^*$ acting on $H^2(\tilde{f}(\widetilde{\Omega}))$ with nonzero eigenvalue $\lambda$.  That is, 
\begin{eqnarray}
\label{pullback-ev}
\left( \tilde{f}^* - \lambda \right)^k u = 0, \quad \mathrm{but} \quad \left( \tilde{f}^* - \lambda \right)^{k-1} u \neq 0
\end{eqnarray}
for some $k \in \ZZ_+$ and $\lambda \neq 0$.  We claim $u$ can be extended to an eigenfunction of $\tilde{f}^*$ acting on $H^2(B)$ with the same eigenvalue.  Indeed, if (\ref{pullback-ev}) holds, we have 
\begin{eqnarray*}
\left( \tilde{f}^* - \lambda \right)^k u = \left[ \sum_{j=0}^k {k \choose j} (-1)^j \lambda^j \left(\tilde{f}^* \right)^{k-j} \right]u = 0,
\end{eqnarray*}
which motivates setting 
\begin{eqnarray}
\label{u-ext}
\tilde{u} := (-1)^{k+1} \lambda^{-k}\left[ \sum_{j=0}^{k-1} {k \choose j} (-1)^{j} \left( \tilde{f}^*\right)^{k-j} \right]u.
\end{eqnarray}
The lowest order pullback on the right hand side of (\ref{u-ext}) is order $1$, and since $\tilde{f}(B) \subset \tilde{f}(\widetilde{\Omega})$, $\tilde{u}$ is in $H^2(B)$.  As $(\tilde{f}^* - \lambda)$ commutes with $\lambda^j (\tilde{f}^*)^{k-j}$, we have 
\begin{eqnarray*}
\left( \tilde{f}^* - \lambda \right)^k \tilde{u} = 0, \quad \mathrm{but} \quad \left( \tilde{f}^* - \lambda \right)^{k-1} \tilde{u} \neq 0.
\end{eqnarray*}
Lastly, if $u$ is a generalized eigenfunction of $\tilde{f}^*$ acting on $H^2(B)$, clearly the restriction of $u$ to $\tilde{f}(\widetilde{\Omega})$ is a generalized eigenfunction with the same eigenvalue for $\tilde{f}^*$ acting on $H^2(\tilde{f}(\widetilde{\Omega}))$.  Then the trace of $\tilde{f}^*$ acting on $H^2(\tilde{f}(\widetilde{\Omega}))$ and $H^2(B)$ are the same and we can apply Lemma \ref{pullback-lemma'} to get (\ref{pullback-formula-statement}).
\end{proof}

\end{document}

%% file: R-fig2.tex
\begin{picture}(0,0)%
\includegraphics{R-fig2.pstex}%
\end{picture}%
\setlength{\unitlength}{3158sp}%
\begingroup\makeatletter\ifx\SetFigFont\undefined%
\gdef\SetFigFont#1#2#3#4#5{%
  \reset@font\fontsize{#1}{#2pt}%
  \fontfamily{#3}\fontseries{#4}\fontshape{#5}%
  \selectfont}%
\fi\endgroup%
\begin{picture}(4599,4074)(1069,-3373)
\put(1351,-3211){\makebox(0,0)[lb]{\smash{{\SetFigFont{10}{12.0}{\rmdefault}{\mddefault}{\updefault}{\color[rgb]{0,0,0}$\reals_{-}$}%
}}}}
\put(4951,-2986){\makebox(0,0)[lb]{\smash{{\SetFigFont{10}{12.0}{\rmdefault}{\mddefault}{\updefault}{\color[rgb]{0,0,0}$i \reals_{+}$}%
}}}}
\put(5176,-1561){\makebox(0,0)[lb]{\smash{{\SetFigFont{10}{12.0}{\rmdefault}{\mddefault}{\updefault}{\color[rgb]{0,0,0}$s=ir$}%
}}}}
\put(2326,-2086){\makebox(0,0)[lb]{\smash{{\SetFigFont{10}{12.0}{\rmdefault}{\mddefault}{\updefault}{\color[rgb]{0,0,0}$R'_{\alpha}$}%
}}}}
\put(3826, 89){\makebox(0,0)[lb]{\smash{{\SetFigFont{10}{12.0}{\rmdefault}{\mddefault}{\updefault}{\color[rgb]{0,0,0}$|\Re s| = |\Im s|^{\alpha}$}%
}}}}
\put(3076,-211){\makebox(0,0)[lb]{\smash{{\SetFigFont{10}{12.0}{\rmdefault}{\mddefault}{\updefault}{\color[rgb]{0,0,0}$D_1(r)$}%
}}}}
\put(4126,-211){\makebox(0,0)[lb]{\smash{{\SetFigFont{10}{12.0}{\rmdefault}{\mddefault}{\updefault}{\color[rgb]{0,0,0}$D_2(r)$}%
}}}}
\put(2926,-1861){\makebox(0,0)[lb]{\smash{{\SetFigFont{10}{12.0}{\rmdefault}{\mddefault}{\updefault}{\color[rgb]{0,0,0}$[-2r^\alpha,0] \times i[r,r+r^\alpha]$}%
}}}}
\put(4426,-2461){\makebox(0,0)[lb]{\smash{{\SetFigFont{10}{12.0}{\rmdefault}{\mddefault}{\updefault}{\color[rgb]{0,0,0}$R_{\alpha}$}%
}}}}
\put(1951,314){\makebox(0,0)[lb]{\smash{{\SetFigFont{10}{12.0}{\rmdefault}{\mddefault}{\updefault}{\color[rgb]{0,0,0}$|\Re s| = 2|\Im s|^{\alpha}$}%
}}}}
\put(1201,-511){\makebox(0,0)[lb]{\smash{{\SetFigFont{10}{12.0}{\rmdefault}{\mddefault}{\updefault}{\color[rgb]{0,0,0}$|\Re s| = 5|\Im s|^{\alpha}$}%
}}}}
\end{picture}%